%% file: osc.tex
\documentclass[10pt,a4paper,reqno]{amsart} 

\usepackage{latexsym,floatflt}
\usepackage{epsfig}
\usepackage{rotating}
\usepackage{amssymb}
\usepackage[T1]{fontenc}
\usepackage{afterpage}
\usepackage{color}

\usepackage{amssymb,amsfonts,amsmath,stmaryrd}

\catcode`\@=11
\def\section{\@startsection{section}{1}%
 \z@{.7\linespacing\@plus\linespacing}{.5\linespacing}%
 {\normalfont\bfseries\scshape\centering}}
\def\subsection{\@startsection{subsection}{2}%
  \z@{.5\linespacing\@plus\linespacing}{.5\linespacing}%
  {\normalfont\bfseries\scshape}}

\def\subsubsection{\@startsection{subsubsection}{3}%
 \z@{.5\linespacing\@plus\linespacing}{-.5em}
  {\normalfont\bfseries\itshape}}
\catcode`\@=12

\newfont{\bbold}{msbm10 scaled \magstep1}
\newfont{\bbolds}{msbm7 scaled \magstep1}

\newcommand{\zss}{\mbox{\bbolds Z}}
\newcommand{\qs}{\mbox{\bbold Q}}

\newcommand{\fps}{formal power series}

\newcommand{\ps}{permutations}

\newcommand{\bx}{\bar x}
\newcommand{\by}{\bar y}
\newcommand{\bu}{\bar u}
\newcommand{\bv}{\bar v}

\newcommand{\GL}{\mathbb{L}}

\newcommand{\GK}{\mathbb{K}}

\newcommand{\Os}{\mathcal O}
\newcommand{\V}{\mathcal V}
\newcommand{\W}{\mathcal W
}\newcommand{\Ref}[1]{(\ref{#1})}
\newcommand{\beq}{\begin{equation}}
\newcommand{\eeq}{\end{equation}}
\newcommand{\gf}{generating function}
\newcommand{\gfs}{generating functions}

\def\emm#1,{{\em #1}}

\def\cqfd{\par\nopagebreak\rightline{\vrule height 3pt width 5pt depth 2pt}
\medbreak}

 \newtheorem{Theorem}{Theorem}
 \newtheorem{Proposition}[Theorem]{Proposition}

\title[Three osculating walkers]
{Three osculating walkers}

\author{Mireille Bousquet-M\'elou}
\address{CNRS, LaBRI, Universit\'e Bordeaux 1, 351 cours de la Lib\'eration,
  33405 Talence Cedex, France}
\email{mireille.bousquet@labri.fr}
\thanks{MBM  was partially supported by the European Commission's IHRP
  Programme, grant HPRN-CT-2001-00272, ``Algebraic Combinatorics in
  Europe''} 
\date{April 8, 2005}

\begin{document}
\maketitle

\begin{flushright}
{To Tony Guttmann, on the occasion of his 60th birthday}
\end{flushright}
\begin{abstract}
We consider three directed walkers on the square lattice, which move
simultaneously at each tick of a clock and never cross. Their
trajectories form a \emm non-crossing configuration of walks,. This
configuration is said to be \emm osculating, if the walkers never
share an edge, and \emm vicious, (or: non-intersecting) if they never
meet. 

We give a closed form expression for the \gf\ of osculating
configurations starting from prescribed points. This \gf\ turns out to
be algebraic. We also relate the enumeration of osculating
configurations with prescribed starting and ending points to the
(better understood) enumeration of non-intersecting configurations.

Our method is based on a step by step decomposition of osculating
configurations, and on the solution of the functional equation
provided by this decomposition. 
   \end{abstract}

\section{Introduction}\label{section-introduction}
\begin{figure}[hbt]
\begin{center}
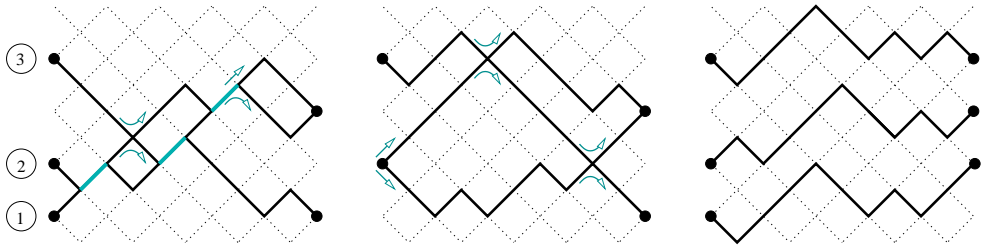
\end{center}
\caption{A non-crossing configuration, an osculating configuration,
  and a non-intersecting (vicious) configuration.}
\label{fig-conf}
\end{figure}

Consider $p$ directed walkers on the (rotated) square lattice,
labelled from 1 to $p$ (Figure~\ref{fig-conf}).  At time 0, all of them
are located at 
abscissa 0, at respective (even) ordinates $j_{0,1}, \ldots , j_{0,p}$, with
$j_{0,1}\le 
j_{0,2} \le \cdots \le j_{0,p}$. Then, at each tick of a clock, each of them
moves to the  
right. More precisely, at each $m \in \llbracket 1, n\rrbracket $,
each walker takes either a 
North-East step $(1,1)$ or a South-East step $(1,-1)$. The set of
trajectories of these walkers, stopped at time $n$, is called a
\emm configuration of paths, of length $n$. This configuration is \emm
non-crossing, if, at each 
time $m$, the ordinates $j_{m,1}, \ldots , j_{m,p}$ of the $p$ walkers
remain ordered as they were at time 0, that is, if $j_{m,1}\le 
j_{m,2} \le \cdots \le j_{m,p}$.  The configuration is \emm
non-intersecting, (or \emm vicious,) if  $j_{m,1}<
j_{m,2} < \cdots < j_{m,p}$ for all $m \in \llbracket 0,n\rrbracket
$. The configuration is \emm osculating, 
if, as soon as $j_{m,k}=j_{m,k+1}$ for some $m \in \llbracket
0,n-1\rrbracket $ and $k \in 
\llbracket 1,p-1\rrbracket $, then  $j_{m+1,k}<j_{m+1,k+1}$. That is,
the walkers are allowed to meet, but cannot share an edge nor cross. In an
osculating configuration, every pair $(m,k)$ such that
$j_{m,k}=j_{m,k+1}$ and $m<n$ is 
called an \emm osculation,. For instance, the second configuration of
Figure~\ref{fig-conf} has 3 osculations (the final contact of the
endpoints is not counted as an osculation). Observe
that in an osculating configuration of positive length, three walkers
never occupy the same site.

Configurations of vicious and osculating walkers have attracted a lot of
attention in the past 20 years, both in combinatorics and in
statistical physics\footnote{Moreover, Tony Guttmann has personally
  observed non-crossing configurations of vacillating and oscillating
  runners in the final portion of the Marathon du M\'edoc.}. Vicious walkers are known to be related to many
important combinatorial objects, like plane
partitions, Young tableaux, symmetric functions, perfect
matchings... to name just a 
few. See~\cite[Ch.~2]{stanley-vol1} and~\cite[Ch.~7]{stanley-vol2}, as
well as~\cite{brenti,gessel-viennot2,stembridge-tspp}, for
instance. In physics, they were 
introduced by Fisher as a model of 
``wetting and melting'', and they can be considered as networks of
polymers~\cite{essam,fisher,melbourne-viennot-no-wall,kratti-one-wall}. The
enumeration of non-intersecting configurations 
of walks is  well understood: in many cases, it reduces to
the evaluation of a determinant~\cite{gessel-viennot}, or a
Pfaffian~\cite{stembridge}, for  
which efficient tools are now available~\cite{kratti-det1}.

Configurations of \emm osculating, walkers naturally arise in physics, in
the ice model (or 6-vertex model)~\cite{baxter}. More recently, it was
realized that 
they are also connected to some famous matrices, called \emm
alternating sign matrices,~\cite{mbm-habsieger}.  These matrices are
renowned for 
having defeated the combinatorial community for more than a decade: it was
conjectured in 1983 that their number is given by a remarkably simple
product formula~\cite{mrr}, but this formula was only proved, with
tremendous 
difficulty, in 1996~\cite{zeilberger-asm}. A bit later, a second
proof was found,   based on some former work on the 6-vertex
model~\cite{kuperberg,izergin}. Let us finally mention that there exists a
conjectural formula for the number of osculating configurations with
fixed endpoints~\cite{brak}. In this formula, the determinant that
usually appears for non-intersecting configurations is replaced by a
more complicated sum on permutations.

\smallskip
In this note, we consider configurations of \emm three,
walkers. Following a terminology inspired by Duplantier~\cite{duplantier}, and
now commonly used  in many physics papers, we say that a non-crossing
configuration
of  three walkers starting respectively at
 ordinates $0, 2i$ and $2(i+j)$, with $i,j \ge 0$, is an
 $(i,j)$\emm-star,. We give explicitly the length \gf\ of \emm osculating,
 $(i,j)$-stars, which 
turns out to be a simple algebraic (quadratic) series
(Proposition~\ref{propo-stars-gf}). 
The case $i=j=1$ of our expression  proves a
conjecture of Guttmann and V\"oge~\cite{guttmann-voge}. We refine our result by
taking into account, in the enumeration, the number of osculations,
thus proving a refined conjecture of Essam~\cite{essam}. We also
obtain the length \gf\ of \emm vicious, $(i,j)$-stars. 

Finally, we consider the enumeration of osculating stars in which the
distances between the three endpoints are also fixed. More
precisely, we study the \gf
$$
\Os _{i,j}(t;x,y) = \sum_{k, \ell, n\ge 0} o _{i,j}^{(k,\ell)}(n)x^k y^\ell
t^n,
$$
where $ o _{i,j}^{(k,\ell)}(n)$ is the number of osculating
$(i,j)$-stars in which the  three walkers  end at time $n$ at
ordinates $j_{1}, 
j_{2}, j_{3}$, with $j_{2}-j_{1}=2k $ and
$j_{3}-j_{2}=2\ell $. We call $\Os _{i,j}$ the 
\emm complete \gf\ , of $(i,j)$-stars. 
 We find an intriguing relation
between this series and the complete \gf\ $\V _{i,j}$ defined
similarly for the (better understood) \emm vicious, walkers
(Proposition~\ref{propo-stars-xy}).  This relation proves that
$\Os(t;x,y)$ is \emm D-finite,, and allows us to compute $ o
_{i,j}^{(k,\ell)}(n)$ explicitly for given values of $i,j, k,
\ell$. In particular, we prove a second conjecture of Guttmann and
V\"oge on the number of osculating \emm watermelons,. 
Note that $\Os _{i,j}(t;1,1)$ is simply the length \gf\
of osculating $(i,j)$-stars. 

\bigskip
Let us conclude this introduction by recalling some definitions and notation
 on formal power series.
 Given a ring $\GL$ and $k$  indeterminates
$x_1, \ldots , x_k$,
 we denote by
 $\GL[x_1, \ldots , x_k]$  
the ring of polynomials in 
$x_1, \ldots , x_k$ 
with coefficients in $\GL$. We denote by  
$\GL[[x_1, \ldots , x_k]]$ 
the ring of formal power series in 
the $x_i$
%
with coefficients in $\GL$.
 A {\em Laurent polynomial\/} in the $x_i$ is a polynomial in both the
 $x_i$ and the $\bx _i=1/x_i$. 
%
For $F\in \GL[[t]]$, we denote by $[t^n]F$ the coefficient of $t^n$ in
$F(t)$. 
 If $F$ is a formal series in $t$ whose coefficients are
Laurent polynomials in $x$, we denote by $F^+$ the {\em positive part of $F$ in
  $x$\/}, that is,
\beq \label{positive-part}
F= \sum_{n\ge 0} t^n \sum_{i\in \zss} f_i(n) x^i \ \Rightarrow \ 
F^+= \sum_{n\ge 0} t^n \sum_{i > 0} f_i(n) x^i.
\eeq
We define similarly  the negative
part of $F$.

Assume, from now on, that $\GL$ is a field. We denote by
 $\GL(x_1, \ldots , x_k)$ the field of rational functions of $x_1,
\ldots , x_k$ with coefficients in $\GL$.
%
A series $F$ in $\GL[[x_1, \ldots , x_k]]$ is {\em algebraic\/}
if there exists a non-trivial polynomial $P$ with coefficients in
$\GL$ such that 
$P(F,x_1, \ldots , x_k)=0.$ 
The sum and product of algebraic series is algebraic.  
%
%
The series $F$   is {\em D-finite\/} 
if the partial derivatives of $F$ span a finite 
dimensional vector space over the field $\GL(x_1, \ldots , x_k)$;
see~\cite{stanleyDF} for the one-variable case,
and~\cite{lipshitz-diag,lipshitz-df} otherwise.  In other 
words, for $1\le 
i\le k$, the series $F$ satisfies a non-trivial partial differential
equation of the form
$$\sum_{\ell=0}^{d_i}P_{\ell,i}\ \frac{\partial ^\ell
F}{\partial x_i^\ell} =0,$$
where $P_{\ell,i}$ is a polynomial in the $x_j$.
Any algebraic series is D-finite.
 The sum and product of D-finite series are D-finite.
 Finally, if $F$ is 
D-finite, then any {\em diagonal\/} of $F$  is also
 D-finite~\cite{lipshitz-diag} (the diagonal of $F$ in $x_1$ and $x_2$ is
 obtained by keeping only those monomials for which the exponents of
 $x_1$ and $x_2$ are equal).  We shall use the following
 consequence of this result: if
 $F(t,x) \in \GL[x,\bx ][[t]]$ is 
algebraic,
then the positive part of
 $F$ in $x$ is D-finite, as well as its negative part.

\section{The complete \gf\ of osculating stars}\label{section-main}
The results stated in this section will be proved in the next
section. Our first proposition  deals with the \emm length, \gf\ of
$(i,j)$-stars. 
\begin{Proposition}
\label{propo-stars-gf}
For $i , j \ge 0$ and $(i,j)\not=(0,0)$, the length \gf \ of
osculating $(i,j)$-stars is algebraic and 
belongs to $\qs(t,\sqrt{1-8t})$.  For instance,
$$
\Os_{1,1}(1,1)= \frac{3-15t-4t^2-3(1-t)\sqrt{1-8t}}{8t^2(1+t)}.
$$
More generally, let $T\equiv T(t)$ be the unique power series in $t$
satisfying $T=2t(1+T)^2$:
$$
T=\frac{1-4t- \sqrt{1-8t}}{4t}.
$$
Then
$$
\begin{array}{lll}
(1-8t)\Os_{i,j}(1,1)&=&\displaystyle
1-3\frac{T^{j+1}}{1+2T}
+3\frac{T^{i+j+1}}{2+T}
-3\frac{T^{i+1}}{1+2T}\\
\\
&= &\displaystyle 1 - 3\, \frac t{1+t} \left( 
T^j(2+T)-T^{i+j}(1+2T) +T^i(2+T)
\right).
\end{array}
$$
For $i , j \ge 0$, the length \gf \ of {\em vicious} $(i,j)$-stars is
algebraic and belongs to $\qs(t,\sqrt{1-8t})$:
$$
(1-8t)\V_{i,j}(1,1)=(1-T^i)(1-T^j).
$$
\end{Proposition}
%
%
These results have also been obtained, independently and via a
different approach, by  Gessel~\cite{ira}. We compare both
approaches after the proof of Proposition~\ref{propo-stars-gf}.
The expression of $\Os_{1,1}(1,1)$ was conjectured in~\cite{guttmann-voge}. 
In Section~\ref{section-refinements}, we refine the above result by
taking into account the 
number of osculations: we prove that the refined \gf \ belongs to
$\qs(t,u,\sqrt{1-8t})$ (where the variable $u$ counts the
osculations). This interpolates between osculating stars and vicious
stars. 

Note that 
$$
1-8t=\frac{(1-T)^2}{(1+T)^2} =\frac{2t} T(1-T)^2.
$$
Hence the above result for vicious $(i,j)$-stars specializes, when
$i=j=1$,  to 
$$
\V_{1,1}(1,1)= \frac T{2t} = \sum_{n\ge 0 } \frac{2^n}{n+2}
  {{2n+2}\choose {n+1}}t^n,$$
as was already proved in~\cite{melbourne-viennot-no-wall}. As
explained there, counting vicious
$(1,1)$ stars is equivalent to counting semi-standard Young tableaux
having at most 3 columns. 

\medskip
For the \emm complete, \gf \ of
stars,
we obtain the following result.
\begin{Proposition}
\label{propo-stars-xy}
For $i , j \ge 0$, the complete \gf \ of osculating $(i,j)$-stars is
D-finite,  and can be expressed in terms of the complete \gfs \ of
vicious stars: 
$$
(1+t) \Os_{i,j}(x,y)= x^iy^j +t \ \frac {x+y+xy}{xy} \Big(
 \V_{i,j}(x,y)+\V_{i+1,j}(x,y)+\V_ {i,j+1}(x,y)\Big).
$$
\end{Proposition}
%
\smallskip
\noindent  Let us make two comments on this result.

\smallskip
\noindent {\bf 1. D-finite series.}
 The number of vicious $(i,j)$-stars of length $n$ such that the
endpoints of the three paths are respectively $-n+2r, -n+2r+2k$ and
$-n+2r+2k+2\ell$ can be expressed, using the Gessel-Viennot
method~\cite{gessel-viennot},  as the following determinant:
$$
v_{i,j}^{(k,\ell)}(r,n) = 
\left|
\begin{array}{ccc}
\displaystyle {n \choose {r}} &\displaystyle {n \choose {r+k }} & \displaystyle {n \choose {r+k +\ell }}\\
\\
\displaystyle {n \choose {r-i}} &\displaystyle {n \choose {r+k -i}} & \displaystyle {n \choose {r+k +\ell -i}}\\
\\
\displaystyle {n \choose {r-i-j}} &\displaystyle {n \choose {r+k -i-j}} &\displaystyle  {n \choose {r+k +\ell -i-j}}
\end{array}
\right|.
$$
Hence the complete \gf\ of vicious $(i,j)$-stars reads
$$
\V_{i,j}(t;x,y)=\sum_{k,\ell, n\ge 0}\sum_{r=0}^ n
v_{i,j}^{(k,\ell)}(r,n) x^k y^\ell
t^n,
$$
and the closure properties of D-finite series~\cite{lipshitz-df} imply that
$\V_{i,j}(t;x,y)$  is D-finite. The expression of
Proposition~\ref{propo-stars-xy} 
shows that $\Os_{i,j}(t;x,y)$ is also D-finite.

\medskip
\noindent
{\bf 2. Watermelons of all sorts.} In particular, when $i=j=1$, we
obtain 
$$
\begin{array}{lllll}
[t^nxy] \V_{1,1}(t;x,y)&=&\displaystyle\sum_{r=0}^ n v_{1,1}^{(1,1)}(r,n)\\
& =&
\displaystyle \frac 2 {(n+1)(n+2)^2}\sum_{r=0}^ n {{n+2} \choose{r}}{{n+2}
  \choose{r+1}}{{n+2} \choose{r+2}}:=b_{n+1}.
\end{array}
$$
The configurations counted by the series $[xy] \V_{1,1}$ are sometimes
called (vicious) \emm watermelons,. The number of watermelons of
length $n$, given 
above, is also the number of Baxter permutations of length $n+1$ 
(see~\cite{dulucq-guibert-baxter} and references therein).
Let us now set  $i=0$ and $j=1$ in
Proposition~\ref{propo-stars-xy}. Since $\V_{0,j}=0$ for all $j$, this gives
$$
(1+t)\Os_{0,1}(x,y)= y+ t\,\frac{x+y+xy}{xy} \,\V_{1,1}(x,y).
$$
Similarly,
$$
(1+t)\Os_{1,0}(x,y)= x+ t\,\frac{x+y+xy}{xy} \,\V_{1,1}(x,y).
$$
Recall that $  \V_{1,1}(x,y)$ is a multiple of $xy$, and
extract  from these two identities the coefficient of $x^1y^0$. This
gives
\beq
\label{1001-baxter}
[x^1y^0] \Os_{0,1}(x,y)= [x^1y^0] \Os_{1,0}(x,y) -\frac 1{1+t}=
\frac t{1+t} [xy]\, \V_{1,1}(t;x,y) =
\frac {B(t)}{1+t},
\eeq
where $B(t) =\sum_{n\ge 1} b_nt^n$ is the \gf\ of Baxter \ps.
Hence 
$$
o_{0,1}^{(1,0)}(n)= \sum_{k= 1} ^n (-1)^{n-k } b_k,
$$
where $b_k$ is the number of Baxter permutations of length $k$. 
Note also that
$$
o_{1,0}^{(1,0)}(n)=o_{0,1}^{(1,0)}(n)+(-1)^n,
$$
which does not seem to be combinatorially obvious.

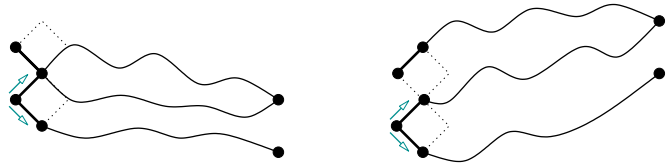
\begin{figure}[hbt]
\begin{center}
\input{watermelons.pstex_t}
\end{center}
\caption{The first three steps in an osculating $(0,1)$-star.}
\label{fig-watermelons}
\smallskip
\end{figure}

Now the first 3-tuple of steps  in an osculating $(0,1)$-star is very
constrained: only two possibilities are allowed for these first steps
(Figure~\ref{fig-watermelons}). This observation implies that 
$$
[x^1y^0] \Os_{0,1}(x,y) =  t[x^1y^0] \Os_{1,0}(x,y)+ t[x^1y^0]
\Os_{1,1}(x,y).$$
From~\Ref{1001-baxter}, we obtain
$$
[x^1y^0] \Os_{1,1}(x,y)= \frac{(1-t)B(t)-t}{t(1+t)}.
$$

\noindent Similarly, there are only two possibilities for the last 3-tuple of
steps in a configuration counted by $[x^1y^0] \Os_{1,1}(x,y)$. This gives:
$$
\begin{array}{llll}
[x^1y^0] \Os_{1,1}(x,y)&= & t[x^0y^1] \Os_{1,1}(x,y)+ t[x^1y^1]
\Os_{1,1}(x,y)\\
&= & t[x^1y^0] \Os_{1,1}(x,y)+ t[x^1y^1] \Os_{1,1}(x,y)& \hbox{by symmetry}.
\end{array}
$$
Hence the \gf\ of ``osculating watermelons'' is finally
\beq
\label{1111-baxter}
[x^1y^1] \Os_{1,1}(x,y)= \frac {1-t}{t^2(1+t)} \Big( (1-t)
B(t)-t\Big).
\eeq
Using the {\sc Maple} packages {\tt EKHAD} and {\tt GFUN}, one can
prove that the 
series $B(t)$ satisfies the following linear differential equation:
$$
12\,t
-6\, \left( 1-2\,t \right)B(t) 
 -2\,t
 \left( 3-14\,t-8\,{t}^{2} \right)B'(t)
 -{t}^{2} \left( t+1 \right)  \left(1- 8\,t \right)
B''(t)=0.
$$
By combining the last two equations, we obtain a differential equation
satisfied by the \gf\ $[x^1y^1]
\Os_{1,1}(t;x,y)$ of osculating watermelons. This equation was
conjectured in~\cite[Eq.~(4.38)]{guttmann-voge}.

\section{Proofs}
\begin{figure}[tb]
\begin{center}
\input{8steps.pstex_t}
\end{center}
\caption{The eight possible moves.}
\label{eight}
\smallskip
\end{figure}
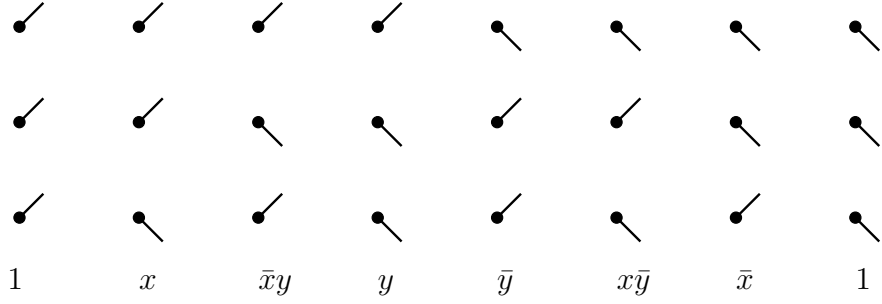
\subsection{Proof of Proposition~\ref{propo-stars-gf}}
For the sake of simplicity, let us denote by $\Os(x,y)$ the complete
\gf \ of osculating $(i,j)$-stars (instead of $\Os_{i,j}(x,y)$). Imagine
we construct these stars step by step as follows: we start from the
star reduced to three points, and add, at each tick of the clock, one
step to each of the three walks. In general, there are $2^3=8$ ways of
adding these steps. The way in which they modify the distances between
the endpoints of the walks is summarized in Figure~\ref{eight}. However, if
two walks end at the same place, exactly $6$ of these $8$ moves are
illegal. This simple construction translates into the following
equation:
$$
 \Os(x,y)=x^iy^j 
+t\left(1+x+\bx y +y +\by +x\by +\bx+1  \right)\Os(x,y)
\hskip 8 cm$$
$$
\hskip 25 mm
-t\left(1+x +\by +x\by +\bx+1\right)\Os(x,0) 
-t\left(1+\bx y +y +\by +\bx+1 \right)\Os(0,y),
$$
which can be rewritten as
$$
\left(xy-t(1+x)(1+y)(x+y)\right) \Os(x,y)
\hskip 10 cm
$$
\vskip -4mm\begin{eqnarray}
\hskip 2cm &
=& x^{i+1}y^{j+1}  -t(x+y+xy)(1+x)\Os(x,0) -t(x+y+xy)(1+y)\Os(0,y),
\nonumber
\\
&
=&x^{i+1}y^{j+1}  -(x+y+xy)P(x) -(x+y+xy)Q(y),
\label{eq-S-ij}
\end{eqnarray}
where $P(x)=t(1+x)\Os(x,0)$ and $ Q(y)=t(1+y)\Os(0,y)$.
We call the coefficient of
$\Os(x,y)$  the {\em kernel\/} $K(x,y)$ of the equation:
\beq
\label{kernel}
K(x,y)=xy-t(1+x)(1+y)(x+y).
\eeq
We are going to apply to~\Ref{eq-S-ij} the {\em
obstinate kernel method\/} that has already been used
in~\cite{bousquet-motifs,bousquet-kreweras}.  
The classical kernel method consists in coupling the
variables $x$ and $y$ so as to cancel the  kernel $K(x,y)$.
This  gives some ``missing'' information about
the series $P(x)$ and $Q(y)$ (see for
instance~\cite{bousquet-petkovsek-recurrences,hexacephale}). In its 
obstinate version, the kernel 
method is combined with a procedure that constructs and exploits
several  (related) couplings $(x,y)$. This procedure is essentially
borrowed from~\cite{fayolle-livre}, where similar functional equations
occur in a probabilistic context.

Let us first fix $x$, and consider the kernel as a quadratic polynomial in
$y$. Its two roots are:
$$
\begin{array}{lclllll}
Y_0(x)&=&\displaystyle 
\frac{1-t(1+x)(1+\bx) -\sqrt{1-2t(1+x)(1+\bx )-t^2(1-x^2)(1-\bx^2)}}
{2t(1+\bx)} \\
&=& \displaystyle (1+x)t+(1+x)^2(1+\bx)t^2 + O(t^3) , \\
\\
Y_1(x)&=&\displaystyle 
\frac{1-t(1+x)(1+\bx) +\sqrt{1-2t(1+x)(1+\bx )-t^2(1-x^2)(1-\bx^2)}}
{2t(1+\bx)} \\
&=& \displaystyle \frac x{1+x}\frac 1 t-(1+x)-(1+x)t+ O(t^2) . \\
\end{array}$$
Observe that $Y_0Y_1=x$.
The first root $Y_0$ is a \fps\ in $t$, and can thus be substituted for $y$
in~\Ref{eq-S-ij}. This gives a functional equation relating 
$P$ and $Q$:
\begin{equation}
P(x)+Q(Y_0)=\frac{x^{i+1}Y_0^{j+1}}{x+Y_0+xY_0}.
\label{kernel0}
\end{equation}
Replacing $y$ by $Y_1$ in $\Os(x,y)$ would not give a  well-defined
power series in $t$, so that we must resist the temptation of this
substitution. However, the following procedure will produce other
interesting pairs $(x,y)$ that cancel the kernel.

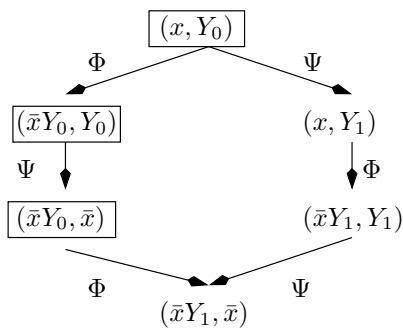
\begin{figure}[hbt]
\begin{center}
\input{diagram1.pstex_t}
\end{center}
\caption{The orbit of  $(x,Y_0)$ under the action of  $\Phi$ and $\Psi$.}
\label{diagram}
\end{figure}

 Let $(X,Y)\not = (0,0)$ be a pair of Laurent series in $t$ with
coefficients in a field $\GK$  such that $K(X,Y)=0$. Recall that $K$
is quadratic in $x$ and $y$. In particular, the equation  $K(x,Y)=0$
admits a second solution $X'$. 
 Define $\Phi(X,Y)= (X',Y)$.
Similarly, define $\Psi(X,Y)= (X,Y')$, where $Y'$ is 
the second solution  of $K(X,y)=0$. Note that $\Phi$ and $\Psi$ are
involutions. Moreover, with the kernel given by~\Ref{kernel}, one has
$Y'=X/Y$ and $X'=Y/X$. Let us examine the action of $\Phi$ and $\Psi$ on the
pair $(x,Y_0)$: we obtain an orbit of cardinality $6$
(Figure~\ref{diagram}).

The 6  pairs of power series given in Figure~\ref{diagram} cancel the
kernel, and we have framed 
the ones that can be legally substituted for $(x,y)$ in the main functional 
equation~\Ref{eq-S-ij}. We thus obtain
 {\em three\/} equations relating the unknown series $P(x)$ and $Q(x)$:
$$
\left\{
\begin{array}{lll}
P(x)+Q(Y_0)&=&\displaystyle
\frac{x^{i+1}Y_0^{j+1}}{x+Y_0+xY_0}, \\
\\
P(\bx Y_0)+Q(Y_0)&=&\displaystyle\frac{\bx^{i}Y_0^{i+j+1}}{1+x+Y_0},\\
\\
P(\bx Y_0)+Q(\bx)&=&\displaystyle\frac{\bx^{i+j}Y_0^{i+1}}{x+Y_0+xY_0}.
\end{array}
\right.
$$
By combining these three equations, we obtain a relation between
$P(x)$ and $Q(\bx)$:
\beq \label{P-Q}
P(x)+Q(\bx)=\frac{x^{i+1}Y_0^{j+1}}{x+Y_0+xY_0}
-\frac{\bx^{i}Y_0^{i+j+1}}{1+x+Y_0}
+\frac{\bx^{i+j}Y_0^{i+1}}{x+Y_0+xY_0}.
\eeq
%
Setting $x=1$ in the above equation gives
$$
P(1)+Q(1)=\frac{T^{j+1}}{1+2T}
-\frac{T^{i+j+1}}{2+T}
+\frac{T^{i+1}}{1+2T},
$$
where $T=Y_0(1)$ is the series defined in Proposition~\ref{propo-stars-gf}.
Setting $x=y=1$ in~\Ref{eq-S-ij} gives
$$
\left(1-8t\right) \Os_{i,j}(1,1)=1-3P(1)-3Q(1).
$$
The first part of Proposition~\ref{propo-stars-gf} follows.

\medskip
We now apply the same approach to the enumeration of vicious
$(i,j)$-stars, or, more precisely, to the enumeration of \emm
quasi-vicious, stars: these are the  osculating $(i,j)$-stars that
are only allowed to meet at their (rightmost) endpoint.
 Let $ \W_{i,j}\equiv \W$
denote their complete \gf . Note that the complete \gf\ of vicious
$(i,j)$-stars is, for $i, j \ge 1$,
$$
\V(x,y)\equiv \V_{i,j}(x,y)=\W_{i,j}(x,y)-\W_{i,j}(x,0)-\W_{i,j}( 0,y).
$$
We can construct quasi-vicious stars  step by step, as we did for
general osculating stars. The difference is that now, no further move
is possible when two walkers are the same place. This gives
$$
 \W(x,y)=x^iy^j +t(1+\bx)(x+y)(1+\by)\left(
\W(x,y)-\W(x,0)-\W(0,y)\right),
$$
that is,
\beq
\label{V-W-eq}
(1-t(1+\bx)(x+y)(1+\by)) \V(x,y) =
x^{i}y^{j} -\W(x,0)-\W(0,y).
\eeq
The rest of the argument copies what we did for osculating stars. In
particular,
\beq\label{W-kernel}
\W(x,0)+\W(0,\bx)=x^{i}Y_0^{j}-\bx^i Y_0^{i+j} +\bx^{i+j}Y_0^i,
\eeq
hence 
$$
\W(1,0)+\W(0,1)=T^{j}-T^{i+j} +T^i,
$$
and the expected expression of $\V(1,1)$ follows using~\eqref{V-W-eq}.
\cqfd

\noindent{\bf Note.} Proposition~\ref{propo-stars-gf} has also been
obtained by Gessel~\cite{ira}. Here, we  sketch his approach and
compare it to ours. Gessel considers the \gf
$$
G(t;u,v)\equiv G(u,v)= \sum_{i,j \ge 0} u^i v^j \Os _{i,j}(t;1,1).
$$
This series counts all stars, by their length and by the position of
their starting points. He then writes a recurrence relation for the
coefficients of $G$, which is equivalent to the following functional
equation:
$$
G(u,v)=\frac 1 {(1-u)(1-v)} -1 + t(1+\bu)(1+\bv)(u+v) G(u,v)
\hskip 40mm$$
$$\hskip 60mm
-t(1+u)(1+\bu+\bv) G(u,0)-t(1+v)(1+\bu+\bv)G(0,v) .
$$
This equation reflects a recursive description of stars based on the
deletion of the \emm first, step of each path. Then, he conjectures
that $G(u,v)$ is a rational function of $u, v$ and the series $T$,
\emm guesses, this rational function with the help of {\sc Maple}, and
finally checks that it satisfies the  functional equation (or
the corresponding recurrence relation on the coefficients of $G$). 

The main difference between his approach and ours is that we 
derive the solution of the functional equation without having to
guess anything. This allows us to generalize easily
Proposition~\ref{propo-stars-gf} in various ways, as shown by
Propositions~\ref{propo-stars-xy} and \ref{propo-stars-gf-osc}.

\subsection{Proof of Proposition~\ref{propo-stars-xy}}
We now wish to evaluate the complete \gf\ of $(i,j)$-stars, not only
their length \gf. Let us go back
to~\Ref{P-Q}. The series $P(x)=t(1+x)\Os(x,0)$ is a \fps \ in $t$ with
coefficients in 
$x\qs[x]$, while $Q(\bx)$ is a \fps \ in $t$ with coefficients in
$\bx\qs[\bx]$. Hence  $P(x)$ and $Q(\bx)$ are respectively  the
positive part and the negative part  of the
right-hand side of~\Ref{P-Q}, as defined by~\eqref{positive-part}. But this
right-hand side  is an algebraic series,
and this implies $P(x)$ and $Q(x)$ are D-finite.  Going back to the
main equation~\Ref{eq-S-ij}, we conclude 
that the complete \gf \ $\Os(x,y)$  is D-finite too.

A similar treatment may be applied to quasi-vicious stars: since
$\W(x,0)$ and $\W(0,\bx)$ are  power series in $t$ with coefficients in
$x\qs[x]$ and $\bx\qs[\bx]$ respectively, it follows from~\Ref{W-kernel} that
they are, respectively, the positive and the negative part of 
$x^{i}Y_0^{j}-\bx^i Y_0^{i+j} +\bx^{i+j}Y_0^i$.

 The only information that we have used to determine
$\W(x,0)$ and $\W(0,x)$ is the fact that,
 for each pair $(X,Y)$ framed in the diagram of
Figure~\ref{diagram},
\beq\label{W-X-Y}
\W_{i,j}(X,0)+\W_{i,j}(0,Y)= X^i Y^j.
\eeq
Similarly, our determination  of $P(x)$ and $Q(y)$ is based on the fact that,
for each such pair $(X,Y)$,
\beq\label{P-Q-X-Y}
P(X)+Q(Y)=\frac{X ^{i+1} Y^{j+1}}{X+Y+XY}.
\eeq
Now observe that, for any pair $(X,Y)$ such that $K(X,Y)=0$,
$$
\frac 1{X+Y+XY}=\frac t{1+t} \frac{1+X+Y}{XY}.
$$
In particular,  the identity~\Ref{P-Q-X-Y} can be rewritten
$$
P(X)+Q(Y)=\frac {tX ^{i} Y^{j}}{1+t} (1+X+Y).
$$
Comparing with~\Ref{W-X-Y} gives, by linearity,
\begin{eqnarray*}
(1+t)P(x)/t&=& \W_{i,j}(x,0) +\W_{i+1,j}(x,0) +\W_{i,j+1}(x,0), \\
(1+t)Q(y)/t&=& \W_{i,j}(0,y) +\W_{i+1,j}(0,y) +\W_{i,j+1}(0,y). 
\end{eqnarray*}
We now plug these expressions of 
$P(x)$ and $Q(y)$ 
into~\Ref{eq-S-ij}, use~\Ref{V-W-eq}, and obtain
$$
(1+t)\Os_{i,j}(x,y)=x^iy^j +t \ \frac {x+y+xy}{xy} \Big(
\V_{i,j}(x,y)+\V_{i+1,j}(x,y)+\V_{i,j+1}(x,y)\Big)
$$ 
as stated in Proposition~\ref{propo-stars-xy}.

\section{The number of osculations}\label{section-refinements}

In this section, we refine the \gf\ of osculating $(i,j)$-stars by adding a
new indeterminate $u$, which keeps track of the number of
osculations. We denote by $\Os_{i,j}(t; u,x,y)\equiv \Os_{i,j}(u,x,y)$
the refined \gf . 
For instance, the $(0,1)$-star of Figure~\ref{fig-conf} has a
contribution $t^{10}x^2y^0u^3$ is this \gf.

\begin{Proposition}
\label{propo-stars-gf-osc}
For $i , j \ge 0$ and $(i,j)\not = (0,0)$, the  \gf \ that counts of
$(i,j)$-stars by their 
length and number of osculations is algebraic and
belongs to $\qs(t,u,\sqrt{1-8t})$.  
%
More precisely, let $T\equiv T(t)$ be the unique power series in $t$
satisfying $T=2t(1+T)^2$. Then
$$
(1-8t)\Os_  {i,j}(1,1)=1-\frac{4-u}{(1+T)^2-uT^2}
\left({T^{j+1}}
-\frac{T^{i+j+1}(2(1+T)-u)}{2(1+T)-uT}
+{T^{i+1}}\right).
$$
\end{Proposition}
\noindent {\bf Proof.} As in the proof of Proposition~\ref{propo-stars-gf}, we
first write a 
functional equation defining $\Os_{i,j}(u,x,y)\equiv \Os(x,y)$. We
have to weight each osculation by $u$, which gives:
\begin{eqnarray*}
 \Os(x,y)&=&x^iy^j 
+t(1+\bx)(1+\by)(x+y)\Os(x,y)\\
&&-t(\bx+\by+1)(1+x)\Os(x,0) 
-t(\bx+\by+1)(1+y)\Os(0,y)\\
&&+ ty(u-1)(1+\bx)\Os(x,0)+tx(u-1)(1+\by)\Os(0,y) .
\end{eqnarray*}
This  can be rewritten as
\beq
K(x,y)\Os(x,y)= x^{i+1}y^{j+1}  -(x+y+xy+y^2(1-u))P(x)
-(x+y+xy+x^2(1-u))Q(y), \label{eq-S-ij-0-u}
\eeq
where $P(x)=t(1+x)\Os(x,0)$ and $ Q(y)=t(1+y)\Os(0,y)$, and the kernel
$K(x,y)$ is still given by~\Ref{kernel}.
The rest of the proof follows the same principles as the proof of
Proposition~\ref{propo-stars-gf}. We successively replace the pair $(x,y)$ by the
three framed pairs of Figure~\ref{diagram}. This gives three linear equations
that relate $P(x), P(\bx Y_0), Q(Y_0)$ and $Q(\bx)$. We eliminate $
P(\bx Y_0)$ and $ Q(Y_0)$ between these three equations to obtain a
linear relationship between $P(x)$ and $Q(\bx)$. Setting $x=1$ and
$Y_0=T$ in this equation gives:
$$
\Big(P ( 1 ) +Q ( 1) \Big)  
\left((1 +T)^2-uT^2 \right) = {T}^{j+1}-
{\frac { \left(2-u+ 2
T \right) {T}^{i+j+1}}{2(1+T)-uT}}+{T}^{i+1}.
$$
But according to~\Ref{eq-S-ij-0-u}, 
$$
(1-8t)\Os(1,1)= 1-(4-u) \left(P ( 1 ) +Q ( 1) \right)  .
$$
Proposition~\ref{propo-stars-gf-osc} follows. The case $i=j=1$ of our
result was conjectured in~\cite{essam}.

\section{Discussion}

The main question raised by this paper is whether this approach can be
extended to more than three walkers. Of course, the step by step
construction can still be exploited: in general, it gives a functional
equation that defines the  \gf\ $\Os_{i_1, \ldots , i_{p-1}}(t;x_1,
\ldots , x_{p-1})$  counting osculating $(i_1, \ldots , i_{p-1})$-stars
according to their length and the distances between their
endpoints. But the problem is how 
to \emm solve, this equation... 

The connection between osculating and vicious walkers stated in
Proposition~\ref{propo-stars-xy} is intriguing. Since we are at a loss to extend it
to more walkers, let us, very modestly, state the corresponding
results for two walkers. 
We use notations that are similar to those used for three walkers, and
should be self-explanatory. We take into account the number of
osculations using an indeterminate $u$. The functional equation that
defines the complete 
\gf\ $\Os_i(t;x)\equiv \Os(x)$ of osculating $i$-stars is
$$
\Os_i(x)= x^i + t(2+x+\bx) \Os_i(x) -t(2+\bx+x(1-u)) \Os_i(0),
$$
that is,
\beq
\label{osc2}
(1-t(1+x)(1+\bx)) \Os_i(x)= x^i  -t(2+\bx+x(1-u)) \Os_i(0).
\eeq
The equation satisfied by the \gf\ $\W_i(t;x)$ of quasi-vicious $i$-stars reads
\beq
\label{vicious2}
(1-t(1+x)(1+\bx)) \V_i(x)= x^i - \W_i(0),
\eeq
where $\V_i(x)= \W_i(x)- \W_i(0)$ is the length \gf \ of vicious
$i$-stars. The standard kernel method gives
$$
\W_i(0)= X^i, \quad \Os_i(0)= \frac {X^{i+1}}{t((1+X)^2-uX^2))}
=\frac {X^{i}}{1-tuX} ,
$$
where $X\equiv X(t)$ is the only power series in $t$ that cancels the
kernel $K(x)=1-t(1+x)(1+\bx)$:
$$
X= \frac{1-2t-\sqrt{1-4t}}{2t}.
$$
Setting $x=1$ in \Ref{osc2} and \Ref{vicious2} gives the counterpart
of Propositions~\ref{propo-stars-gf} and~\ref{propo-stars-gf-osc}:
$$
(1-4t) \Os_i(t;1)=  1 -  \frac {(4-u)X^{i+1}}{(1+X)^2-uX^2}
= 1 -  \frac {(4-u)tX^{i}}{1-tuX}\quad \hbox{ and }
\quad 
(1-4t) \V_i(t;1)=  1- X^i.
$$
To obtain a relation between \emm complete, \gfs \ of osculating and
vicious stars, we observe that
$$
(1-tuX)(1-u+2tu+tuX)= (1+ut)^2-u.
$$
Hence the above expression of $\Os_i(0)$ can be rewritten as
$$
\Os_i(0)= \frac {X^i(1-u+2tu+tuX)}{(1+ut)^2-u}= 
\frac 1 {(1+ut)^2-u} \Big((1-u+2tu) \W_i(0) +
tu\W_{i+1}(0)\Big) .
$$
We now plug this expression back in~\Ref{osc2} and use~\Ref{vicious2}
to obtain the following counterpart of Proposition~\ref{propo-stars-xy}:
\begin{multline*}
  ((1+ut)^2-u) \Os_i(x)=\\
 x^i(1-u+2tu+xtu(1-u)) + t\,\frac {1+2x+x^2(1-u)} x \Big((1-u+2tu) \V_i(x) +
tu\V_{i+1}(x)\Big).
\end{multline*}
When $u=1$, this specializes to 
$$
(2+t) \Os_i(x)= 2x^i + t\,\frac {1+2x} x \Big(2 \V_i(x) +
\V_{i+1}(x)\Big).
$$

Finally, it would  be interesting to find purely combinatorial
proofs of Propositions~\ref{propo-stars-gf} and~\ref{propo-stars-xy}.
\bigskip
\noindent
%



\bibliographystyle{plain}
\bibliography{biblio}


\end{document}

%% file: osc.pstex_t
\begin{picture}(0,0)%
\includegraphics{osc.pstex}%
\end{picture}%
\setlength{\unitlength}{2171sp}%
\begingroup\makeatletter\ifx\SetFigFont\undefined%
\gdef\SetFigFont#1#2#3#4#5{%
  \reset@font\fontsize{#1}{#2pt}%
  \fontfamily{#3}\fontseries{#4}\fontshape{#5}%
  \selectfont}%
\fi\endgroup%
\begin{picture}(11134,2766)(125,-2494)
\end{picture}%

%% file: watermelons.pstex_t
\begin{picture}(0,0)%
\includegraphics{watermelons.pstex}%
\end{picture}%
\setlength{\unitlength}{2171sp}%
\begingroup\makeatletter\ifx\SetFigFont\undefined%
\gdef\SetFigFont#1#2#3#4#5{%
  \reset@font\fontsize{#1}{#2pt}%
  \fontfamily{#3}\fontseries{#4}\fontshape{#5}%
  \selectfont}%
\fi\endgroup%
\begin{picture}(7515,1838)(729,-1314)
\end{picture}%

%% file: 8steps.pstex_t
\begin{picture}(0,0)%
\includegraphics{8steps.pstex}%
\end{picture}%
\setlength{\unitlength}{3947sp}%
\begingroup\makeatletter\ifx\SetFigFont\undefined%
\gdef\SetFigFont#1#2#3#4#5{%
  \reset@font\fontsize{#1}{#2pt}%
  \fontfamily{#3}\fontseries{#4}\fontshape{#5}%
  \selectfont}%
\fi\endgroup%
\begin{picture}(5497,1880)(526,-3569)
\put(526,-3511){\makebox(0,0)[lb]{\smash{\SetFigFont{12}{14.4}{\familydefault}{\mddefault}{\updefault}{\color[rgb]{0,0,0}$1$}%
}}}
\put(1351,-3511){\makebox(0,0)[lb]{\smash{\SetFigFont{12}{14.4}{\familydefault}{\mddefault}{\updefault}{\color[rgb]{0,0,0}$x$}%
}}}
\put(2851,-3511){\makebox(0,0)[lb]{\smash{\SetFigFont{12}{14.4}{\familydefault}{\mddefault}{\updefault}{\color[rgb]{0,0,0}$y$}%
}}}
\put(3601,-3511){\makebox(0,0)[lb]{\smash{\SetFigFont{12}{14.4}{\familydefault}{\mddefault}{\updefault}{\color[rgb]{0,0,0}$\by$}%
}}}
\put(4351,-3511){\makebox(0,0)[lb]{\smash{\SetFigFont{12}{14.4}{\familydefault}{\mddefault}{\updefault}{\color[rgb]{0,0,0}$x\by$}%
}}}
\put(5101,-3511){\makebox(0,0)[lb]{\smash{\SetFigFont{12}{14.4}{\familydefault}{\mddefault}{\updefault}{\color[rgb]{0,0,0}$\bx$}%
}}}
\put(5851,-3511){\makebox(0,0)[lb]{\smash{\SetFigFont{12}{14.4}{\familydefault}{\mddefault}{\updefault}{\color[rgb]{0,0,0}$1$}%
}}}
\put(2101,-3511){\makebox(0,0)[lb]{\smash{\SetFigFont{12}{14.4}{\familydefault}{\mddefault}{\updefault}{\color[rgb]{0,0,0}$\bx y$}%
}}}
\end{picture}

%% file: diagram1.pstex_t
\begin{picture}(0,0)%
\includegraphics{diagram1.pstex}%
\end{picture}%
\setlength{\unitlength}{3947sp}%
\begingroup\makeatletter\ifx\SetFigFont\undefined%
\gdef\SetFigFont#1#2#3#4#5{%
  \reset@font\fontsize{#1}{#2pt}%
  \fontfamily{#3}\fontseries{#4}\fontshape{#5}%
  \selectfont}%
\fi\endgroup%
\begin{picture}(2247,1965)(1129,-1561)
\put(2101,239){\makebox(0,0)[lb]{\smash{\SetFigFont{10}{12.0}{\familydefault}{\mddefault}{\updefault}\special{ps: gsave 0 0 0 setrgbcolor}$(x,Y_0)$\special{ps: grestore}}}}
\put(1651, 14){\makebox(0,0)[lb]{\smash{\SetFigFont{10}{12.0}{\familydefault}{\mddefault}{\updefault}\special{ps: gsave 0 0 0 setrgbcolor}$\Phi$\special{ps: grestore}}}}
\put(3376,-661){\makebox(0,0)[lb]{\smash{\SetFigFont{10}{12.0}{\familydefault}{\mddefault}{\updefault}\special{ps: gsave 0 0 0 setrgbcolor}$\Phi$\special{ps: grestore}}}}
\put(1201,-661){\makebox(0,0)[lb]{\smash{\SetFigFont{10}{12.0}{\familydefault}{\mddefault}{\updefault}\special{ps: gsave 0 0 0 setrgbcolor}$\Psi$\special{ps: grestore}}}}
\put(3001, 14){\makebox(0,0)[lb]{\smash{\SetFigFont{10}{12.0}{\familydefault}{\mddefault}{\updefault}\special{ps: gsave 0 0 0 setrgbcolor}$\Psi$\special{ps: grestore}}}}
\put(1651,-1411){\makebox(0,0)[lb]{\smash{\SetFigFont{10}{12.0}{\familydefault}{\mddefault}{\updefault}\special{ps: gsave 0 0 0 setrgbcolor}$\Phi$\special{ps: grestore}}}}
\put(2926,-1411){\makebox(0,0)[lb]{\smash{\SetFigFont{10}{12.0}{\familydefault}{\mddefault}{\updefault}\special{ps: gsave 0 0 0 setrgbcolor}$\Psi$\special{ps: grestore}}}}
\put(1201,-361){\makebox(0,0)[lb]{\smash{\SetFigFont{10}{12.0}{\familydefault}{\mddefault}{\updefault}\special{ps: gsave 0 0 0 setrgbcolor}$(\bx Y_0,Y_0)$\special{ps: grestore}}}}
\put(1201,-961){\makebox(0,0)[lb]{\smash{\SetFigFont{10}{12.0}{\familydefault}{\mddefault}{\updefault}\special{ps: gsave 0 0 0 setrgbcolor}$(\bx Y_0,\bx)$\special{ps: grestore}}}}
\put(3001,-361){\makebox(0,0)[lb]{\smash{\SetFigFont{10}{12.0}{\familydefault}{\mddefault}{\updefault}\special{ps: gsave 0 0 0 setrgbcolor}$(x,Y_1)$\special{ps: grestore}}}}
\put(3001,-961){\makebox(0,0)[lb]{\smash{\SetFigFont{10}{12.0}{\familydefault}{\mddefault}{\updefault}\special{ps: gsave 0 0 0 setrgbcolor}$(\bx Y_1,Y_1)$\special{ps: grestore}}}}
\put(2101,-1561){\makebox(0,0)[lb]{\smash{\SetFigFont{10}{12.0}{\familydefault}{\mddefault}{\updefault}\special{ps: gsave 0 0 0 setrgbcolor}$(\bx Y_1,\bx)$\special{ps: grestore}}}}
\end{picture}

%% file: osc.bbl
\begin{thebibliography}{10}

\bibitem{hexacephale}
C.~Banderier, M.~Bousquet-M{\'e}lou, A.~Denise, P.~Flajolet, D.~Gardy, and
  D.~Gouyou-Beauchamps.
\newblock Generating functions for generating trees.
\newblock {\em Discrete Math.}, 246(1-3):29--55, 2002.

\bibitem{baxter}
R.~J. Baxter.
\newblock {\em Exactly solved models in statistical mechanics}.
\newblock Academic Press Inc., London, 1989.

\bibitem{bousquet-motifs}
M.~Bousquet-M\'elou.
\newblock Four classes of pattern-avoiding permutations under one roof:
  generating trees with two labels.
\newblock {\em Electronic J. Combinatorics}, 9(2):Research Paper 19, 2003.

\bibitem{bousquet-kreweras}
M.~Bousquet-M{\'e}lou.
\newblock Walks in the quarter plane: Kreweras' algebraic model.
\newblock {\em Ann. Appl. Probab.}, to appear.

\bibitem{mbm-habsieger}
M.~Bousquet-M\'elou and L.~Habsieger.
\newblock {Sur les matrices \`a signes alternants.}
\newblock {\em Discrete Math.}, 139(1-3):57--72, 1995.

\bibitem{bousquet-petkovsek-recurrences}
M.~Bousquet-M{\'e}lou and M.~Petkov{\v{s}}ek.
\newblock Linear recurrences with constant coefficients: the multivariate case.
\newblock {\em Discrete Math.}, 225(1-3):51--75, 2000.

\bibitem{brak}
R.~Brak.
\newblock Osculating lattice paths and alternating sign matrices.
\newblock In {\em Formal Power Series and Algebraic Combinatorics}, Wien,
  Austria, 1997.
\newblock See also the addendum on R. Brak's web page.

\bibitem{brenti}
F.~Brenti.
\newblock Determinants of super-{S}chur functions, lattice paths, and dotted
  plane partitions.
\newblock {\em Adv. Math.}, 98(1):27--64, 1993.

\bibitem{dulucq-guibert-baxter}
S.~Dulucq and O.~Guibert.
\newblock Baxter permutations.
\newblock {\em Discrete Math.}, 180(1-3):143--156, 1998.

\bibitem{duplantier}
B.~Duplantier.
\newblock Statistical mechanics of polymer networks of any topology.
\newblock {\em J. Statist. Phys.}, 54(3-4):581--680, 1989.

\bibitem{essam}
J.~W. Essam.
\newblock Three attractive osculating walkers and a polymer collapse
  transition.
\newblock {\em J. Statist. Phys.}, 110(3--6):1191--1207, 2003.

\bibitem{fayolle-livre}
G.~Fayolle, R.~Iasnogorodski, and V.~Malyshev.
\newblock {\em Random walks in the quarter-plane: Algebraic methods, boundary
  value problems and applications}, volume~40 of {\em Applications of
  Mathematics}.
\newblock Springer-Verlag, Berlin, 1999.

\bibitem{fisher}
M.~E. Fisher.
\newblock {Walks, walls, wetting, and melting.}
\newblock {\em J. Stat. Phys.}, 34:667--730, 1984.

\bibitem{ira}
I.~Gessel.
\newblock Personal communication.

\bibitem{gessel-viennot2}
I.~Gessel and G.~Viennot.
\newblock Determinants, paths, and plane partitions.
\newblock Preprint 1989.

\bibitem{gessel-viennot}
I.~Gessel and G.~Viennot.
\newblock {Binomial determinants, paths, and hook length formulae.}
\newblock {\em Adv. Math.}, 58:300--321, 1985.

\bibitem{melbourne-viennot-no-wall}
A.~J. Guttmann, A.~L. Owczarek, and X.~G. Viennot.
\newblock Vicious walkers and {Y}oung tableaux. {I}. {W}ithout walls.
\newblock {\em J. Phys. A}, 31(40):8123--8135, 1998.

\bibitem{guttmann-voge}
A.~J. Guttmann and M.~V{\"o}ge.
\newblock Lattice paths: vicious walkers and friendly walkers.
\newblock {\em J. Statist. Plann. Inference}, 101(1-2):107--131, 2002.

\bibitem{izergin}
A.~G. Izergin.
\newblock {Partition function of a six-vertex model in a finite volume.}
\newblock {\em Sov. Phys. Dokl.}, 32(11):878--879, 1987.

\bibitem{kratti-det1}
C.~Krattenthaler.
\newblock {Advanced determinant calculus.}
\newblock {\em S\'emin. Lothar. Comb.}, 42, 1999.
\newblock B42q, 67 p.

\bibitem{kratti-one-wall}
C.~Krattenthaler, A.~J. Guttmann, and X.~G. Viennot.
\newblock Vicious walkers, friendly walkers and {Y}oung tableaux. {II}. {W}ith
  a wall.
\newblock {\em J. Phys. A}, 33(48):8835--8866, 2000.

\bibitem{kuperberg}
G.~Kuperberg.
\newblock {Another proof of the alternating-sign matrix conjecture.}
\newblock {\em Int. Math. Res. Not.}, 1996(3):139--150, 1996.

\bibitem{lipshitz-diag}
L.~Lipshitz.
\newblock The diagonal of a {$D$}-finite power series is {$D$}-finite.
\newblock {\em J. Algebra}, 113(2):373--378, 1988.

\bibitem{lipshitz-df}
L.~Lipshitz.
\newblock D-finite power series.
\newblock {\em J. Algebra}, 122:353--373, 1989.

\bibitem{mrr}
W.~H. Mills, D.~P. Robbins, and H.~Rumsey.
\newblock {Alternating sign matrices and descending plane partitions.}
\newblock {\em J. Comb. Theory, Ser. A}, 34:340--359, 1983.

\bibitem{stanleyDF}
R.~P. Stanley.
\newblock Differentiably finite power series.
\newblock {\em European J. Combin.}, 1:175--188, 1980.

\bibitem{stanley-vol1}
R.~P. Stanley.
\newblock {\em Enumerative combinatorics. {V}ol. 1}, volume~49 of {\em
  Cambridge Studies in Advanced Mathematics}.
\newblock Cambridge University Press, Cambridge, 1997.
\newblock With a foreword by Gian-Carlo Rota, Corrected reprint of the 1986
  original.

\bibitem{stanley-vol2}
R.~P. Stanley.
\newblock {\em Enumerative combinatorics $2$}, volume~62 of {\em Cambridge
  Studies in Advanced Mathematics}.
\newblock Cambridge University Press, Cambridge, 1999.

\bibitem{stembridge}
J.~R. Stembridge.
\newblock {Nonintersecting paths, pfaffians, and plane partitions.}
\newblock {\em Adv. Math.}, 83(1):96--113, 1990.

\bibitem{stembridge-tspp}
J.~R. Stembridge.
\newblock The enumeration of totally symmetric plane partitions.
\newblock {\em Adv. Math.}, 111(2):227--243, 1995.

\bibitem{zeilberger-asm}
D.~Zeilberger.
\newblock {Proof of the alternating sign matrix conjecture.}
\newblock {\em Electron. J. Comb.}, 3(2), 1996.
\newblock Research paper R13, 84 p.

\end{thebibliography}
